\documentclass[12pt]{amsart}
\usepackage{amsfonts, amsmath, amsthm}
\usepackage{epsfig}
\usepackage{psfrag}

\newtheorem{pro}{Proposition}[section]
\newtheorem{thm}[pro]{Theorem}
\newtheorem{lem}[pro]{Lemma}

\newtheorem{cor}[pro]{Corollary}

\theoremstyle{definition}
\newtheorem{dfn}[pro]{Definition}

\theoremstyle{remark}

\newsavebox{\savepar}

\title{Non-parallel essential surfaces in knot complements} 
\date{\today}
\address{Mathematics Department, California Polytechnic State University}
\email{dbachman@calpoly.edu}
\author{David Bachman}
\begin{document}
\begin{abstract}
We show that if a knot or link has $n$ thin levels when put in thin position then its exterior contains a collection of $n$ disjoint, non-parallel, planar, meridional, essential surfaces. A corollary is that there are at least $n/3$ tetrahedra in any triangulation of the complement of such a knot. 
\end{abstract}
\maketitle


\section{Introduction}
In 1987 D. Gabai introduced the concept of {\it thin position} for knots \cite{gabai:87} to solve the Property $R$ conjecture. Since then it has been used to solve many important questions in 3-manifold topology \cite{gl:89}, \cite{st:93}, \cite{thompson:94}. More recently thin position has become an object of study in itself \cite{thompson:97}, \cite{hk:97}, \cite{sr:02}, \cite{wu:03}.

{\it Thin levels} are particular spheres which appear in a thin presentation of a knot or link. In \cite{thompson:97} A. Thompson shows that if the exterior of a knot or link $K$ does not contain any planar, meridional, essential surfaces then thin position for $K$ has no thin levels. The idea of the proof is to show that the property of being {\it non-trivial}, possessed by a thin level, cannot disappear after any number of compressions. Therefore, compressing any thin level as much as possible yields a planar, meridional, essential surface. 

In this paper we show that the property of being {\it non-parallel}, possessed by a collection of thin levels, cannot disappear after any number of compressions (as long as the compressions are chosen in a suitably nice way). Hence, compressing all thin levels as much as possible yields at least as many non-parallel essential surfaces as the number of thin levels. As a corollary we show that there are at least $n/3$ tetrahedra in any triangulation of the complement of a knot or link with $n$ thin levels. 

The author would like to thank Saul Schleimer for many helpful conversations.

\section{Thin Position}

Suppose $K \subset S^3$ is an arbitrary knot or link with no trivial components and $h$ is some standard height function on $S^3$ (so that for each $p \in (0,1)$, $h^{-1}(p)$ is a 2-sphere), which is a Morse function when restricted to $K$. Let $\{ q'_j \}$ denote the critical values of $h$ restricted to $K$, let $q_j$ be some point in the interval $(q'_j, q'_{j+1})$, and let $S_j=h^{-1}(q_j)$. The following terminology is standard in {\it thin position} arguments (see \cite{gabai:87}).

The {\it width} of $K$ is defined to be the quantity $\sum \limits _j |K \cap S_j|$. A knot is said to be in {\it thin position} if $h$ is chosen so that the width of $K$ is minimal (see \cite{gabai:87}). If $j$ is such that $|K \cap S_j|<|K \cap S_{j-1}|$ and $|K \cap S_j|<|K \cap S_{j+1}|$ then we say the surface $S_j$ is a {\it thin level} of $K$. In other words, a thin level is one which appears just above a maximum of $K$ and just below a minimum.

\begin{dfn}
Suppose $K$ is a knot or link in $S^3$, $h$ is the standard height function, and $\gamma$ is a 1-manifold in the exterior of $K$. We say $\gamma$ is {\it horizontal (with respect to $K$)} if it is contained in a thin level of $K$. If $\gamma$ has endpoints on distinct levels of $h$ then we say it is {\it vertical} if its interior has no critical points (with respect to $h$). If $\gamma$ has endpoints on the same level of $h$ then we say it is {\it $U$-shaped} if its interior has exactly one maximum or minimum. Finally, $\gamma$ is {\it simple} if it is vertical or $U$-shaped.
\end{dfn}

\begin{lem}
\label{l:h-v}
Suppose $K_1$ and $K_2$ are isotopic knots or links in $S^3$ which agree on some 1-manifold, $\alpha$. $K_1$ is not in thin position if any of the following hold:
\begin{enumerate}
    \item $K_2 \setminus \alpha$ is horizontal with respect to $K_1$.
    \item $K_2 \setminus \alpha$ is vertical, but $K_1 \setminus \alpha$ is not.
    \item $K_2 \setminus \alpha$ is $U$-shaped, but $K_1 \setminus \alpha$ is not, and the minimum (maximum) of $K_2 \setminus \alpha$ is at least as high (low) as the minimum (maximum) of $K_1 \setminus \alpha$.
\end{enumerate}
\end{lem}

The horizontal case is proved in \cite{thinball}. The proofs in the vertical and $U$-shaped cases are similar.

\section{Heegaard Splittings}

A {\it compression body} $W$ is a 3-manifold which can be obtained by starting with some surface $F$ (not necessarily closed or connected), forming the product $F \times I$, attaching some number of 2-handles to $F \times \{1\}$, and capping off all remaining 2-sphere boundary components with 3-balls. The surface $F \times \{0\}$ is referred to as $\partial _+ W$. The surface $\partial _- W$ is defined by the equation $\partial W=\partial _+W \cup (\partial F \times I) \cup \partial _- W$. A compression body is {\it non-trivial} if it is not a product. 

A surface $F$ in a 3-manifold $M$ is a {\it Heegaard splitting of M} if $F$ separates $M$ into two compression bodies, $W$ and $W'$, such that $F=\partial _+W=\partial _+ W'$. 
We define $\partial _- M=\partial _-W \cup \partial _- W'$.

We now present the crucial example of a Heegaard splitting for our purposes. Let $T_1$ and $T_2$ denote consecutive thin levels of a knot or link $K$. Let $N$ be the submanifold of $S^3$ cobounded by $T_1$ and $T_2$. Let $S$ be a level 2-sphere which separates the maxima of $K$ in $N$ from the minima. Then $S \backslash K$ is a Heegaard splitting of $N \backslash K$ and $\partial _- (N \backslash K)=T_1 \backslash K \cup T_2 \backslash K$.

For the proof of Lemma \ref{l:GoodMaxCompSeq} will will need the following result:

\begin{lem}
\label{l:haken}
{\rm (Haken \cite{haken:68})} Let $S$ be a Heegaard splitting of a 3-manifold $N$. If $D$ is a compressing disk for $\partial _- N$ then there is a compressing disk $E$ such that $\partial E=\partial D$ and $E \cap S$ is a loop.
\end{lem}

\section{Compressing sequences}

In this section we prove a few preliminary lemmas concerning sequences of compressions of the thin levels of a knot or link. 

\begin{dfn}
Suppose $K$ is a knot or link in $S^3$. A sphere in $S^3$ is {\it trivial} if it is disjoint from $K$, or bounds a ball which contains a single, unknotted arc of $K$. 
\end{dfn}

\begin{dfn}
Suppose $K$ is a knot or link in $S^3$ and $\Sigma$ denotes the union of the thin levels of $K$. A {\it compressing sequence} for $K$ is a sequence of surfaces $\{\Sigma _i\}_{i=0}^n$ such that $\Sigma _0=\Sigma$ and $\Sigma _i$ is obtained from $\Sigma _{i-1}$ by a compression in the complement of $K$. 
\end{dfn}

\begin{dfn}
The compressing sequence $\{\Sigma _i\}_{i=0}^n$ is {\it maximal} if $\Sigma _n$ is incompressible in the complement of $K$.
\end{dfn}

\begin{lem}
\label{l:thompson}
Suppose $K$ is a knot or link in thin position. If $\{\Sigma _i\}_{i=0}^n$ is a maximal compressing sequence for $K$ then each component of $\Sigma _n$ is {\it essential}, {\it i.e.} incompressible and non-trivial.
\end{lem}

\begin{proof}
First we show that each component of $\Sigma _n$ is incompressible in the complement of $K$. Suppose $D$ is a compressing disk for some such component. Then by a standard innermost disk argument we can isotope $D$ so that it intersects every component of $\Sigma_n \backslash K$ in essential loops. Let $\alpha$ denote such a loop which is innermost on $D$. Then $\alpha$ bounds a subdisk of $D$ which is a compressing disk for $\Sigma _n \backslash K$. This contradicts the maximality of $\{\Sigma _i\}_{i=0}^n$.

The remainder of the argument is essentially one of Thompson's from \cite{thompson:94}. We recall this here. The only remaining possibility is that all components of $\Sigma _n$ are trivial. Let $S$ be an innermost component of $\Sigma _n$. Let $D$ be a disk in $S^3$ such that $\partial D=\delta \cup \gamma$, where $D \cap K=\delta$ and $D \cap S=\gamma$. 

Note that $S$ is the result of compressing some thin level $T$ of $K$ some number of times. It follows that $S \cap T$ is connected and we may assume that $\gamma \subset S \cap T$. Hence, we can use $D$ to guide an isotopy of $\delta$ which, in the end, is horizontal with respect to $K$. This contradicts Lemma \ref{l:h-v}.
\end{proof}

\begin{lem}
\label{l:nonUshaped}
Suppose $K$ is a knot or link in thin position. Let $\{\Sigma _i\}_{i=0}^n$ be a compressing sequence for $K$. Let $N_0$ denote the closure of some component of $S^3 \backslash \Sigma _0$. If some component $K'$ of $K \cap N_0$ is $U$-shaped then for all $i$ both points of $\partial K'$ lie on the same component of $\Sigma _i$.
\end{lem}

\begin{proof}
Without loss of generality, assume that $K'$ has a minimum. Then we may isotope $K'$, preserving the width of $K$, so that any other minimum on any other component of $K \cap N_0$ is below the minimum of $K'$. 

Now, suppose the lemma is false. Let $m$ denote the largest integer such that both points of $\partial K'$ lie on the same component of $\Sigma _m$. Then $\Sigma _{m+1}$ is obtained from $\Sigma _m$ by compressing along a disk $D$ whose boundary separates the boundary points of $K'$. Let $N$ denote the closure of the component of $S^3 \backslash \Sigma _m$ which contains $K'$. So $\partial D \subset \partial N$. Note that the interior of $D$ must lie outside of $N$. Otherwise, $D$ would have to intersect $K'$, as $\partial D$ separates $\partial K'$. Let $B$ be either of the balls bounded by $D \cup \partial N$ whose interior lies outside of $N$.

As the component $F$ of $\partial N$ that contains $\partial K'$ is obtained by compressing a thin level, $T$, it must be the case that $F \cap T$ is connected. Hence, we may choose a horizontal arc $\gamma$ in $F \cap T$ which connects the endpoints of $K'$. Furthermore, we may choose such an arc so that it meets $D$ in precisely one point of $\partial D$.

We now perform the isotopy depicted in Figure \ref{f:Kbar}, which can be described as follows:

\begin{enumerate}
    \item Shrink $B$ to a small ball $B'$ at the end of $K'$.
    \item Contract $K'$, pulling $B'$ along with it.
    \item Push $B'$ along the arc $\gamma$.
    \item Inflate $B'$ back to $B$.
\end{enumerate}

        \begin{figure}[htbp]
        \psfrag{k}{$\gamma$}
        \psfrag{K}{$K'$}
        \psfrag{T}{$T$}
        \psfrag{S}{$S$}
        \psfrag{B}{$B$}
        \psfrag{D}{$D$}
        \vspace{0 in}
        \begin{center}
        \epsfxsize=5 in
        \epsfbox{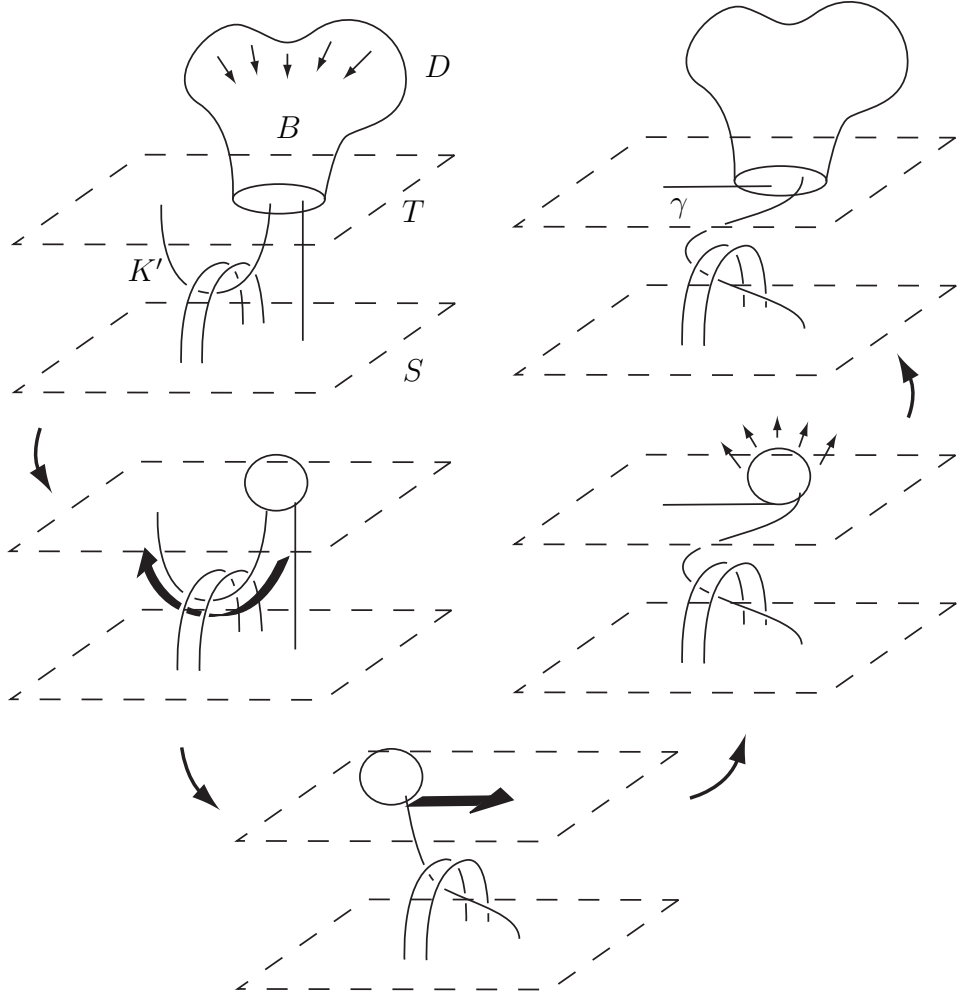}
        \caption{Replacing $K'$ with $\gamma$.}
        \label{f:Kbar}
        \end{center}
        \end{figure}

As in Figure \ref{f:Kbar}, this isotopy will affect the other arcs of $K \cap N_0$ which meet $B$. Let $S$ denote a level 2-sphere which is just below the minimum of $K'$. Note that Steps 2 and 3 of the isotopy take place entirely in the region between $T$ and $S$. Recall that the minimum of $K'$ is above all other minima of the components of $K \cap N_0$. Hence, the subarcs of $K$ that are between $T$ and $S$ and which meet $B$ are vertical. One such arc is depicted in Figure \ref{f:Kbar}. Note that after the isotopy this arc is still vertical. Hence, as in Case 1 of Lemma \ref{l:h-v}, we have reached a contradiction because we have made a non-horizontal subarc of $K$ into a horizontal one while preserving the width everywhere else. 
\end{proof}

\begin{dfn}
The compressing sequence $\{\Sigma _i\}_{i=0}^n$ is {\it good} if for each $i$, each component $N$ of $S^3-\Sigma_i$, and each component $K'$ of $K \cap {\overline N}$ there is a simple arc in ${\overline N}$ connecting the endpoints of $K'$.
\end{dfn} 

\begin{lem}
\label{l:GoodMaxCompSeq}
Suppose $K$ is a knot or link in $S^3$. Then there exists a good maximal compressing sequence for $K$.
\end{lem}

\begin{proof}
To establish the lemma we will define a much more rigid compressing sequence called a {\it Haken sequence} and prove that every Haken sequence is good. We will then show that a Haken sequence of maximal length is a maximal compressing sequence. 

We say a compressing sequence $\{\Sigma _i\}_{i=0}^n$ is {\it Haken} if for each $i$ each component $N$ of $S^3 \backslash \Sigma_i$ contains a 2-sphere $S$ such that 
\begin{enumerate}
    \item $S \backslash K$ is a Heegaard splitting of $N \backslash K$. 
    \item The surface $S$ is obtained from a level surface of $h$ by some sequence of compressions.
    \item Every component of $(K \cap N) \backslash S$ is simple.
\end{enumerate}    
Note that a compressing sequence with a single element, by definition consisting of the thin levels of $K$, is a Haken sequence. Between any two thin levels there is a level $S$ which separates the maxima of $K$ from the minima which has the desired properties. 

We now show that every Haken sequence is good. Let $\{\Sigma _i\}_{i=0}^n$ denote a Haken sequence. Choose some $i$ and let $N$ denote the closure of a component of $S^3-\Sigma_i$. Let $K'$ denote a component of $K \cap N$. To show that our sequence is good we must produce a simple arc in $N$ connecting the endpoints of $K'$. Since our sequence is Haken there is a 2-sphere $S$ in $N$ such that $S \backslash K$ is a Heegaard splitting for $N \backslash K$.

Let $K_1$ and $K_2$ denote the closure of the subarcs of $K' \backslash S$ which contain the points of $\partial K'$. Since $S$ was obtained from a level surface of $h$ by some sequence of compressions we may choose a horizontal arc $\alpha$ in $S$ which connects $K_1 \cap S$ to $K_2 \cap S$ (see Figure \ref{f:gamma}).

        \begin{figure}[htbp]
        \psfrag{a}{$\partial N$}
        \psfrag{b}{$S$}
        \psfrag{c}{$\partial N$}
        \psfrag{d}{$\gamma$}
        \psfrag{1}{$K_2$}
        \psfrag{2}{$\alpha$}
        \psfrag{3}{$K _1$}
        \vspace{0 in}
        \begin{center}
        \epsfxsize=5 in
        \epsfbox{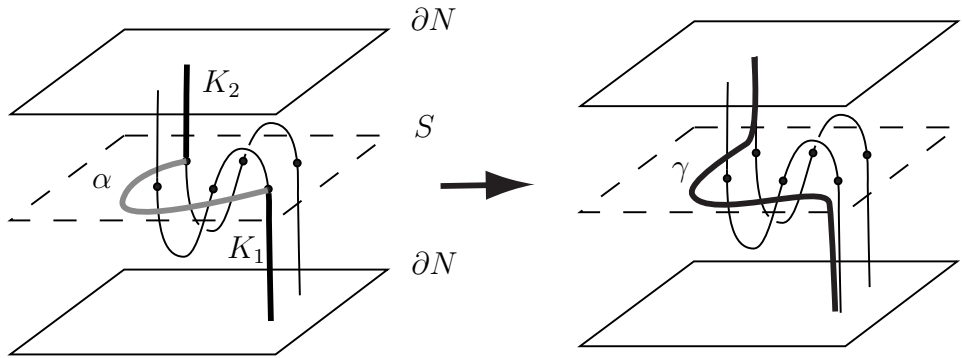}
        \caption{Constructing a simple arc $\gamma$ in $N$.}
        \label{f:gamma}
        \end{center}
        \end{figure}

Since $\{\Sigma _i\}$ is a Haken sequence the arcs $K_1$ and $K_2$ must be vertical. Since the arc $\alpha$ is horizontal we may perturb the arc $\gamma=K_1 \cup \alpha \cup K _2$ to be simple (again, see Figure \ref{f:gamma}). 

We have now shown that that every Haken sequence is good. What remains is to show is that there is a maximal compressing sequence which is a Haken sequence. Again note that for every knot or link $K$ there exists at least one Haken sequence, namely the sequence with one element consisting of the thin levels of $K$. 

We now assume that $\{\Sigma _i\}_{i=0}^n$ is a Haken sequence of maximal length such that some component $T$ of $\Sigma _n$ is compressible in the complement of $K$ ({\it i.e.} $\{\Sigma _i\}_{i=0}^n$ is not a maximal compressing sequence). 
Let $D$ be a compressing disk for $T$ in the complement of $K$. By an innermost disk argument, we may assume that all loops of $D \cap \Sigma _n$ are essential on $\Sigma _n \backslash K$. Let $D'$ be the subdisk of $D$ bounded by an innermost such loop. Let $N$ denote the closure of the component of $S^3 \backslash \Sigma _n$ which contains $D'$. Then $\partial D'$ is a compressing disk for $\partial N$, in the complement of $K$. 

As $\{\Sigma _i\}_{i=0}^n$ is a Haken sequence there is a 2-sphere $S$ in $N$ such that $S \backslash K$ is a Heegaard splitting for $N \backslash K$. By Lemma \ref{l:haken} there is a compressing disk $E$ for $\partial N$, in the complement of $K$, such that $\partial E=\partial D'$ and $E \cap S$ is a simple closed curve, $\delta$. We now compress $S$ along the subdisk of $E$ bounded by $\delta$ to obtain the spheres $S'$ and $S''$ and compress $\Sigma _n$ along $E$ to obtain $\Sigma _{n+1}$ (see Figure \ref{f:compress}). Note that $\{\Sigma _i\}_{i=0}^{n+1}$ is also a Haken sequence, contradicting the maximality of the length of our original choice.
  
        \begin{figure}[htbp]
        \psfrag{D}{$E$}
        \psfrag{2}{$S$}
        \psfrag{N}{$\partial N$}
        \vspace{0 in}
        \begin{center}
        \epsfxsize=5 in
        \epsfbox{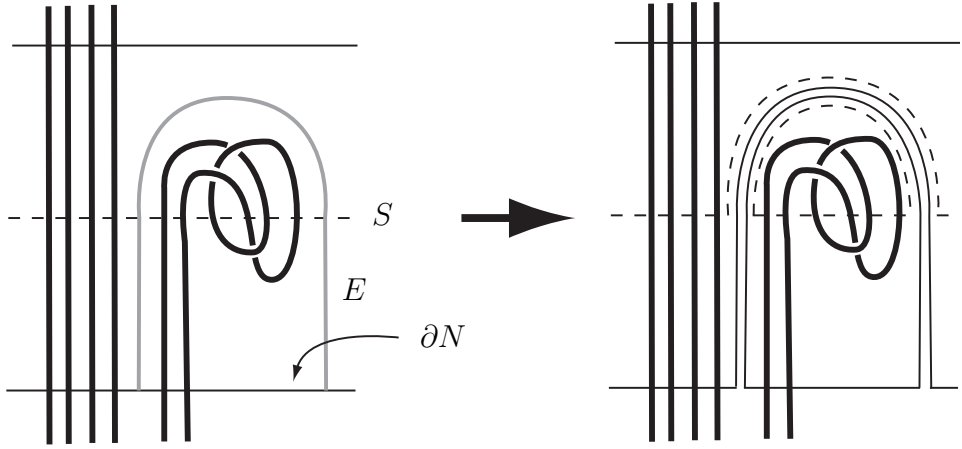}
        \caption{Compressing $S$ and $\Sigma _n$.}
        \label{f:compress}
        \end{center}
        \end{figure}

\end{proof}

\begin{lem}
\label{l:vert}
Suppose $K$ is a knot or link in thin position. Let $\{\Sigma _i\}_{i=0}^n$ be a good compressing sequence for $K$. Let $N$ be the closure of some component of $S^3 \backslash \Sigma_n$ and let $S$ be a component of $\partial N$. Let $K'$ denote the components of $K \cap N$ which meet $S$. If all of the components of $K'$ are parallel and connect distinct components of $\partial N$ then each component of $K'$ is vertical.
\end{lem}

\begin{proof}
To prove the lemma we show that the components of $K'$ must be simple. There are then two possibilities: either they are $U$-shaped or vertical. The former is ruled out by Lemma \ref{l:nonUshaped} and the latter is the desired conclusion.

By way of contradiction, assume the components of $K'$ are not simple. As $\{\Sigma _i\}_{i=0}^n$ is good we may choose a collection of parallel simple arcs $\overline K'$ in $N$ such that for each component $\alpha' \subset K'$ there is a component $\overline \alpha ' \subset \overline K'$ with $\partial \overline \alpha '=\partial \alpha '$.  We now show that $K$ is isotopic to a knot or link $\overline K$ which contains $\overline K'$, such that $K \backslash K'=\overline K \backslash \overline K'$. This then contradicts Lemma \ref{l:h-v}. 

Let $B$ be the ball bounded by $S$ on the side opposite $N$. Note that only one endpoint of each arc component of $K'$ meets $B$, since each such arc connects distinct components of $\partial N$. The isotopy is illustrated in Figure \ref{f:Kbar2}, in the case where $\partial K'$ is contained in a single thin level $T$. The steps are the same in the case where the components of $K'$ connect distinct thin levels. They are as follows:

\begin{enumerate}
    \item Shrink $B$ to a small ball $B'$ at the end of $K'$.
    \item Contract the arcs of $K'$, pulling $B'$ along with it.
    \item Push $B'$ along the arcs of $\overline K'$.
    \item Inflate $B'$ back to $B$.
\end{enumerate}

        \begin{figure}[htbp]
        \psfrag{k}{$\overline K'$}
        \psfrag{K}{$K'$}
        \psfrag{T}{$T$}
        \psfrag{B}{$B$}
        \vspace{0 in}
        \begin{center}
        \epsfxsize=5 in
        \epsfbox{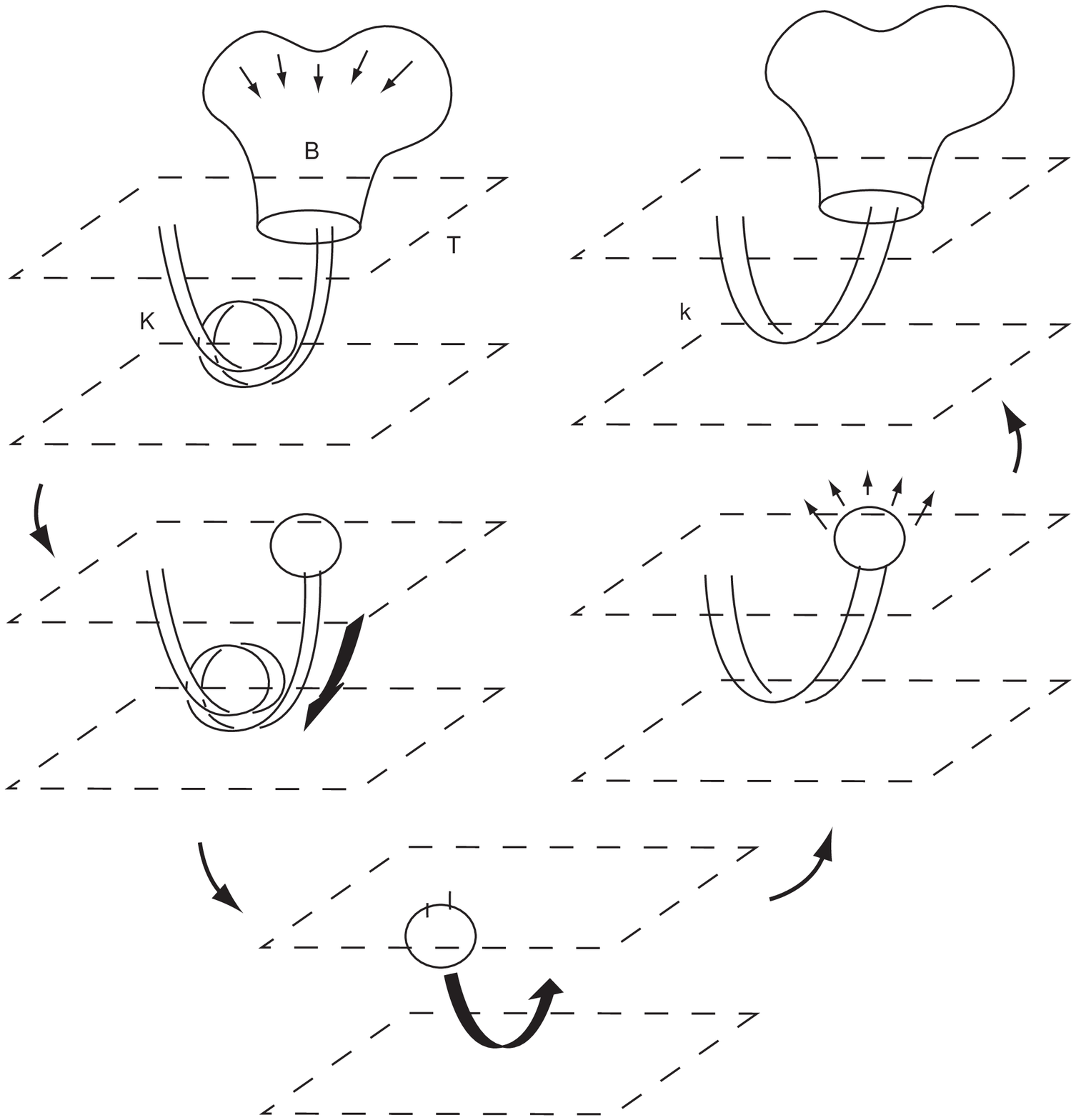}
        \caption{Replacing $K'$ with $\overline K'$.}
        \label{f:Kbar2}
        \end{center}
        \end{figure}

Note that in the case that $\overline K'$ is $U$-shaped we may do a further width-preserving isotopy to make the minima (maxima) of $\overline K'$ appear above (below) the minima (maxima) of $K'$. This is necessary to appeal to Case 3 of Lemma \ref{l:h-v}.
\end{proof}

\section{The Main Theorem.}
In this section we prove our main theorem. 

\begin{thm}
\label{t:main}
If a knot or link has $n$ thin levels when put in thin position then its exterior contains a collection of $n$ disjoint, non-parallel, planar, meridional, essential surfaces.
\end{thm}

\begin{proof}
Let $K$ be a knot or link in thin position with $n$ thin levels. By Lemma \ref{l:GoodMaxCompSeq} we may choose a good maximal compressing sequence $\{\Sigma _i\}_{i=0}^n$ for $K$.

Lemma \ref{l:thompson} implies that the elements of $\Sigma _n$ are essential in the complement of $K$. Let $\mathcal S$ denote a collection of spheres in $S^3$ such that
\begin{enumerate}
    \item every element of $\Sigma _n$ is parallel, in the complement of $K$, to an element of $\mathcal S$ and
    \item no two elements of $\mathcal S$ are parallel in the complement of $K$. 
\end{enumerate}

Our goal is to show that $\mathcal S$ has at least $n$ elements. Let $\Gamma$ denote the dual graph of $\mathcal S$ in $S^3$. $\Gamma$ is then a tree, whose edges correspond to elements of $\mathcal S$. As the number of vertices minus the number of edges of any tree is 1, it suffices to show that $\Gamma$ has at least $n+1$ vertices. 

Note that it is implicit in the assumption that $K$ is in thin position that we have fixed a height function, $h$, on $S^3$. Let $x$ and $y$ denote two critical points of $K$ with respect to $h$, which are separated by a thin level of $K$. Then $x$ and $y$ are separated by an element, $T$, of $\Sigma _n$. We claim that $x$ and $y$ are also separated by an element of $\mathcal S$. 

Suppose this is not the case and let $M$ denote the component of $S^3 \setminus \mathcal S$ which contains $x$ and $y$, so that $T \subset M$. As $T$ is not an element of $\mathcal S$ it must be parallel (in the complement of $K$) to some element of $\mathcal S$, and hence, to some component, $S$, of $\partial M$. Either $x$ or $y$ lies between $T$ and $S$. Assume the former. Since $T$ and $S$ are parallel in the complement of $K$, and $x$ is a point of $K$ which lies between them, $x$ must lie on a subarc, $\alpha$, of $K$ which connects $T$ to $S$. But Lemma \ref{l:vert} implies that $\alpha$ is vertical, contradicting the fact that it contains the critical point, $x$. 

We conclude by noting that our assumption that $K$ had $n$ thin levels implies that there is a collection of $n+1$ critical points of $K$ such that any two are separated by a thin level. The above argument then shows that each of these points must lie in a distinct component of $S^3 \setminus \mathcal S$, implying that $\Gamma$ has at least $n+1$ vertices. 
\end{proof}

\begin{cor}
\label{t:tetra}
Let $K$ be a knot or link which has a thin presentation with $n$ thin levels. Let $t$ be the smallest number of tetrahedra necessary to triangulate the complement of $K$. Then $t \ge \frac{n}{3}$.
\end{cor}

\begin{proof}
In \cite{finite} we give an improvement over the classical Kneser-Haken Finiteness Theorem \cite{kneser:29}, \cite{haken:68} and show that in closed manifolds the size of any collection of pairwise disjoint, closed, essential, 2-sided surfaces is at most twice the number of tetrahedra, $|T|$. Although we do not explicitly state a result there for manifolds (and surfaces) with boundary, the same proof shows that if $M$ is a 3-manifold with non-empty boundary and $\mathcal S$ is a collection of properly embedded, pairwise disjoint, 2-sided, incompressible and boundary incompressible surfaces then $2|\mathcal S| \le g +6|T|$, where $g$ is the maximum number of twisted $I$-bundles that can disjointly embed in $M$. In the complement of the knot or link $K$ in $S^3$ we have $g=0$, so $2|\mathcal S| \le 6t$, or $|\mathcal S| \le 3t$. Now, Theorem \ref{t:main} says there exists such a collection $\mathcal S$ such that $n \le |\mathcal S|$. Hence, $n \le 3t$.
\end{proof}

\bibliographystyle{alpha}

\end{document}